\documentclass{article}
\usepackage{graphicx}
\usepackage{amsmath}
\usepackage{longtable}
\usepackage{program}

\author{Hal Finkel \\ \scriptsize{hfinkel@anl.gov}}
\title{An Iterated, Multipoint Differential Transform Method for Numerically Evolving PDE IVPs}

\begin{document}
\maketitle
\begin{abstract}
Traditional numerical techniques for solving time-dependent partial-differential-equation (PDE) initial-value problems (IVPs) store a truncated representation of the function values and some number of their time derivatives at each time step. Although redundant in the $dx \rightarrow 0$ limit, what if spatial derivatives were also stored? This paper presents an iterated, multipoint differential transform method (IMDTM) for numerically evolving PDE IVPs. Using this scheme, it is demonstrated that stored spatial derivatives can be propagated in an efficient and self-consistent manner; and can effectively contribute to the evolution procedure in a way which can confer several advantages, including aiding solution verification. Lastly, in order to efficiently implement the IMDTM scheme, a generalized finite-difference stencil formula is derived which can take advantage of multiple higher-order spatial derivatives when computing even-higher-order derivatives. As is demonstrated, the performance of these techniques compares favorably to other explicit evolution schemes in terms of speed, memory footprint and accuracy.
\end{abstract}

\section{Introduction}
In traditional methods for solving time-dependent partial-differential-equation (PDE) initial-value problems (IVPs), a truncated representation of the function values and some number of their time derivatives, on some number of spatial slices, are stored on a discrete grid. As is well known, the grid may consist of regularly-spaced points or may be an irregular point collection, and there are many choices for the basis used to represent the function values at the grid points. Nevertheless, evolving the PDE inevitably requires reconstructing spatial derivative values, and because these values are not stored explicitly, they must be reconstructed from the truncated representation of the function values~\footnote{For spectral methods that use a Fourier basis, this reconstruction is trivial. Unfortunately, spectral methods are often impractical for problems with nonlinearities, irregular geometries or nonperiodic boundary conditions.}. This paper reports on an investigation into the construction of a scheme which, unlike traditional schemes, stores some number of spatial derivatives on the grid along with the function values. Although redundant in the $dx \rightarrow 0$ limit, at least in the context of the scheme presented here, storing spatial derivative values will be shown to have practical advantages\footnote{There has been limited investigation of using stored spatial derivatives in the context of finite-volume schemes for evolving systems of conservation laws. In the finite-volume literature, these are known as multi-moment schemes (see~\cite{Ii2007} and references therein).}. To construct a scheme that can make use of, and evolve, stored spatial derivatives, it is natural to consider methods based on generating local power-series expansions of the solution.

Constructing power-series solutions to ordinary differential equations is a basic technique included in most introductory texts on mathematical methods or differential-equations (for example, see~\cite{Boyce2000}). Algorithmically constructing power-series solutions to ordinary differential equations has been extensively studied (see~\cite{Corliss1982}, and references therein), but has remained a little-known technique. Although sometimes used for the numerical evaluation of special functions~\cite{Gil2007}, the construction of power-series solutions has been generally thought of as an analytic tool and not as the basis for numerical algorithms. This is changing, and algorithmically constructing power-series solutions to ordinary differential equations is gaining in popularity. This is now often called the Differential Transform(ation) Method (DTM)~\cite{Chen1996}, but is also known as the Parker-Sochacki method~\cite{Parker1996}, or generically as the Picard-iteration method or the Taylor-series method. In this paper, I adopt the notations and terminology from the DTM literature. The conventions which define the DTM can be straightforwardly extended to multivariate power series for application to PDEs, but applying the resulting definitions to solve PDEs numerically requires dealing with several complexities\footnote{In the context of finite-volume schemes, using local (multivariate) power-series expansions to enhance the accuracy of Riemann solvers forms the basis of the well-known Arbitrary accuracy DERivative (ADER) Riemann technique~\cite{Titarev2002}. Unlike ADER schemes, the technique presented in this paper can be applied to any kind of time-dependent PDE.}. This paper details a practical method for numerically applying the DTM to time-dependent PDEs: an iterated, multipoint DTM (IMDTM), and shows that the resulting evolution algorithm possesses several highly-desirable properties. Notably, by storing a sufficiently-large set of spatial derivatives, the resulting system can be configured to possess a directly-calculable self-consistency constraint, the violation of which measures the quality of the solution.

The remainder of this paper is organized as follows: Section~\ref{sect:ehod} reviews how the evolution equations for higher-order derivatives can be derived. Section~\ref{sect:dtm} reviews the multivariate DTM formalism. The IMDTM algorithm is presented in Section~\ref{sect:imdtm}. The efficient implementation of IMDTM requires an efficient polynomial interpolation scheme capable of using higher-order derivative information. Such a scheme is also presented in Section~\ref{sect:imdtm}. The paper then concludes, and a set of appendices follow.

\section{Evolution Equations for Higher-Order Derivatives}
\label{sect:ehod}
Before discussing the conventions which define the DTM, consider the following two examples which explicitly demonstrate how the higher-order derivatives of a PDE's solution obey PDEs derivable from the PDEs for the solution function. First, consider the 2-D wave equation:
\begin{equation}
	\frac{\partial^2f}{\partial t^2} - \frac{\partial^2 f}{\partial x^2} = 0
	\label{eq:2dwave}
\end{equation}
from which the evolution equation for the higher-order spatial derivatives can be derived by applying the operator $\frac{d^n}{dx^n}$:
\begin{equation}
	\frac{\partial^2}{\partial t^2} \frac{d^nf}{dx^n} - \frac{\partial^2}{\partial x^2} \frac{d^nf}{dx^n} = 0 \, .
	\label{eq:2dwaveddx}
\end{equation}
Because the wave equation is linear, all of the derivatives obey structurally-identical equations. For nonlinear equations, this is not the case. For example, consider the modified KdV equation:
\begin{equation}
	\frac{\partial f}{\partial t} + f^2 \frac{\partial f}{\partial x} + \frac{\partial^3 f}{\partial x^3} = 0
	\label{eq:nlpde}
\end{equation}
for which the evolution equation for the higher-order spatial derivatives are:
\begin{equation}
	\frac{\partial}{\partial t} \frac{d^nf}{dx^n} + \sum_{k=0}^n {n \choose k} \left ( \sum_{j=0}^k {k \choose j} \frac{d^jf}{dx^j} \frac{d^{k-j}f}{dx^{k-j}} \right ) \frac{\partial^{n-k+1} f}{\partial x^{n-k+1}} + \frac{\partial^3}{\partial x^3} \frac{d^nf}{dx^n} = 0 \, .
	\label{eq:nlpdeddx}
\end{equation}
To avoid a ``combinatorial explosion'' in the number of terms, I have used (multivariate) versions of the Leibniz rule (the product rule) and the Fa\`a di Bruno formula (the generalized chain rule)~\cite{Hardy2006}. The multivariate DTM encapsulates these formulas in easy-to-apply recursion relations, thus avoiding the direct use of combinatorial calculus formulas as was done in Equation~\ref{eq:nlpdeddx}. Because the higher-order coefficients depend on the lower-order coefficients, the higher-order coefficients will tend to be slower to compute, and in some cases less numerically stable, than the lower-order coefficients.

It is possible to evolve all of the higher-order coefficients using a traditional PDE evolution scheme by discretizing the associated evolution equations, which will, generally, increase the computational time required by at least a factor of $n$ for evolving $n$ spatial derivatives. The IMDTM scheme allows multiple coefficients to be evolved together, significantly reducing the computational overhead.

\section{Multivariate Differential Transform Method}
\label{sect:dtm}
In many cases, the application of the classic integral transform methods (i.e. Laplace and Fourier) can be reduced to use of a table of substitutions~\cite{Boyce2000}. Similarly, a differential equation can be transformed into a recursion relation for the coefficients of its power series solution(s) using a table of substitutions~\cite{Corliss1982}, and recent literature refers to this as The Differential Transform Method~\cite{Chen1996}. The conventions which define the DTM can be extended to the multivariate case in the obvious way~\cite{Kurnaz2005}, specifically:
\begin{equation}
	W(k_1, k_2, \ldots, k_n) = \frac{1}{k_1!k_2! \cdots k_n!} \left [ \frac{\partial^{k_1 + k_2 + \cdots + k_n}}{\partial x_1^{k_1} \partial x_2^{k_2} \cdots \partial x_n^{k_n}} w(x_1, x_2, \ldots, x_n) \right ]_{0, 0, \ldots, 0}
\label{eq:mdtm}
\end{equation}
so the inverse transform is:
\begin{equation}
	w(w_1, w_2, \ldots, w_n) = \sum_{k_1=0}^\infty \sum_{k_2=0}^\infty \cdots \sum_{k_n=0}^\infty W(k_1, k_2, \ldots, k_n) \prod_{i=1}^n x_i^{k_i}
\label{eq:invmdtm}
\end{equation}
and the basic operations are given in Table~\ref{tab:basemdtm}. By combining Equations~\ref{eq:mdtm} and~\ref{eq:invmdtm}, it is clear that the multivariate DTM is nothing more than the construction of a multivariate Taylor expansion. The efficient evaluation of nonlinear terms is an important practical concern. Fortunately, efficiently computing nonlinear functions of power series is a well-researched problem by the creators of automatic-differentiation software, and is generally done either by directly evaluating a set of multivariate recurrence relations~\cite{Tsukanov2000} or by combining different evaluations of univariate recurrence relations using different directional projections of the original series~\cite{Griewank2000}~\cite{Neidinger2004}. Direct evaluations of the multivariate recurrence relations will be used here, and the DTM translations for additional (nonlinear) functions are provided in Appendix~\ref{sect:mdtmrec}.

\begin{table}[h]
\begin{center}
\begin{tabular}{|l|p{0.6\textwidth}|}
\hline
\bf{Original Function} & \bf{Transformed Function} \\
\hline
$w(x) = y(x) \pm z(x)$ & $W(k) = Y(k) \pm Z(k)$ \\
$w(x) = \lambda y(x)$, $\lambda$ a constant   & $W(k) = \lambda Y(k)$ \\
$w(x) = \frac{\partial^{r_1 + \cdots + r_n}}{\partial x_1^{r_1} \cdots \partial x_n^{r_n}} y(x)$ & $W(k) = \frac{(k+r)!}{k!} Y(k + r)$ \\
$w(x) = y(x) z(x)$ & $W(k) = \sum_{l_1=0}^{k_1} \cdots \sum_{l_n=0}^{k_n} Y(l) Z(k-l)$ \\
$w(x) = x_1^{m_1} \cdots x_n^{m_n}$ & $W(k) = \prod_{i=1}^n \delta_{k_i,m_i}$, $\delta$ is the Kronecker delta \\
\hline
\end{tabular}
\end{center}
\caption{Basic operations under Multidimensional DTM}
\label{tab:basemdtm}
\end{table}

Note that in the example 1+1-dimensional\footnote{1+1-dimensional means that the system has one spatial dimension and one time dimension.} systems discussed in later sections, the temporal coordinate is called $t$ and the spatial coordinate is called $x$. Correspondingly, $k$ is the DTM-coefficient index for temporal derivatives, and $h$ is the DTM-coefficient index for spatial derivatives. Also, an overdot will be used to indicate a time derivative.

\section{Iterated, Multipoint DTM}
\label{sect:imdtm}
Although the multivariate DTM has been applied as a convenient way to analytically construct a power-series solution to PDEs~\cite{Kurnaz2005}~\cite{Jafari2010}, a novel contribution of this paper is a method for using the multivariate DTM as a numerical evolution scheme which can handle an arbitrarily-large Cauchy surface. When solving ODEs, the initial conditions can be exactly represented using a finite set of values. For PDEs, however, the initial conditions are themselves (differentiable) functions and they cannot be represented, in general, using only a finite set of numbers. As expected, applying the multivariate DTM to a PDE results in a recurrence relation which allows higher-order derivative coefficients in one variable to be calculated in terms of the complete set of derivatives with respect to the remaining variables (i.e. on the Cauchy surface). For example, consider the 2-D wave equation:
\begin{equation}
	\frac{\partial^2f}{\partial t^2} = \frac{\partial^2 f}{\partial x^2}
\end{equation}
for which the DTM gives the recurrence:
\begin{equation}
	F(k,h) = \frac{(h+2)(h+1)}{k(k-1)} F(k-2,h+2) \, .
	\label{eq:wave-rec}
\end{equation}
Given the initial tower of derivative coefficients, $F(0, h)$ and $F(1, h)$ (the spatial derivatives of $f$ and $\dot{f}$, respectively), any $F(k, h)$ can be calculated.
 
Even if the radius of convergence of a power series is infinite, when working at finite precision, the useful range may be limited, even as more terms are used. As a result, the IMDTM uses power-series expansions around multiple points on the Cauchy surface to precisely represent the initial conditions. In other words, at each point on the spatial hypersurface, a tower of DTM coefficients are stored, $\{F(k, h)_{t,x}\}$. For a PDE which is first-order in time, only the $k=0$ coefficients are stored; for an equation which is second-order in time, only the $k \in (0,1)$ coefficients are stored; and so on for higher-order PDEs. In combination with the recursion relation derived from the PDE, the power-series expansion of the solution at each spatial point can be reconstructed from the stored coefficients. This provides the IMDTM scheme with a powerful feature: it can yield an internal self-consistency condition. Specifically, the power series around some point must, by self-consistency, provide the values of the function at the neighboring points. As the system is evolved forward in time (via iteration), this condition will break down, and this breakdown can be used to measure the quality of the evolved solution.

As you might expect from the definition of the Taylor expansion, a time step from $t = 0$ to $t = dt$ is calculated as:
\begin{align}
	F(0, 0)_{t = dt, x} & = {} & f(x, dt) & = {} & \sum_{k=0}^\infty F(k, 0)_{t = 0,x} dt^k \label{eq:tsprop00} \\
	F(1, 0)_{t = dt, x} & = {} & \dot{f}(x, dt) & = {} & \sum_{k=0}^\infty (k+1) F(k+1, 0)_{t = 0,x} dt^k \label{eq:tsprop10} \\
	F(0, 1)_{t = dt, x} & = {} & \frac{d}{dx} f(x, dt) & = {} & \sum_{k=1}^\infty F(k, 1)_{t = 0,x} dt^k \label{eq:tsprop01} \\
	\cdots \nonumber
\end{align}
where, again, the overdot indicates a time derivative.

At any time, $t$, the system is specified by an infinite set of coefficients, thus the description must be truncated for numerical implementation. This truncation immediately generates an impediment to the construction of an iterative scheme: at whatever order the truncation is performed, even-higher-order terms are necessary to compute the evolution of the highest-order terms being stored. If these missing coefficients are simply taken to be zero, then the scheme will quickly destabilize. Qualitatively, if $n$ is the highest order stored, and all terms of an order greater than $n$ are taken to be zero, then the term of order $n$ will remain constant, the term of order $n-1$ will grow linearly with time, the term of order $n-2$ will grow quadratically with time, and so on. For the scheme to be practical, a method for computing these missing coefficients must be provided.

\subsection{Interpolating Polynomials}
In many traditional grid-based numerical evolution schemes, the derivatives are approximated by a finite-difference computation. Conceptually, this approximation uses some number of neighboring values to construct an interpolating polynomial, and uses the derivative of that polynomial (at the central point) to estimate the derivative of the function represented on the grid~\cite{Press1997}. Since the IMDTM stores truncated spatial power series on the grid, instead of just the function values, it is possible to use multiple coefficients per point (the lower-order derivative values) in order to reconstruct the needed higher-order derivative values. However, there are several complications which need to be discussed.

First, because the coefficients are stored only to some fixed precision, not all of the neighboring-point power series coefficients should be used for the derivative-reconstruction procedure. To understand this, assume that the power series at a point can be used to calculate the function value at any neighboring point to some precision, and that the stored function value at any neighboring point is accurate to that same precision. This means that at least the neighboring zeroth-order components have no significant information to contribute to the higher-order derivative interpolation. To put it another way, since it takes $n$ orders (i.e. $n$ terms of the power series) to converge to the desired precision, that limited precision is equivalent to introducing an error term of order $n+1$. Adding numbers with an error of order $n+1$ to calculate an order-$(n+1)$ coefficient would be impossible. Even restricting the interpolation procedure to use only the last few orders, it is still, in general, not possible to reconstruct the missing derivatives to the same precision as the inputs. Luckily, by making $dt$ small enough, reconstructing the missing higher-order coefficients to the same precision as the inputs is unnecessary.

The naive method for constructing a polynomial interpolant using higher-order derivative information  -- by solving the implied matrix equation -- is not practical for higher-dimensional systems. Fortunately, it is possible to write down an explicit solution to the general interpolation problem, and use that solution to calculate the required coefficients. The method used here, detailed in Appendix~\ref{sect:ipoly}, uses the multipoint Taylor expansion as defined by L\'opez and Temme~\cite{Lopez2004} to explicitly construct the interpolating polynomial. Although the method seems relatively complicated, a number of the needed factors depend only on the geometry, specifically the relative distances between neighboring grid points, and can be cached once calculated (or, precalculated, if the problem geometry is static). After all of the geometry-dependent factors have been calculated, a linear combination of the input coefficients from the various neighboring points yields the missing higher-order coefficients. This method can be thought of as a generalization of a finite-difference stencil that can make use of higher-order derivative information.

\begin{figure}[h!tb]
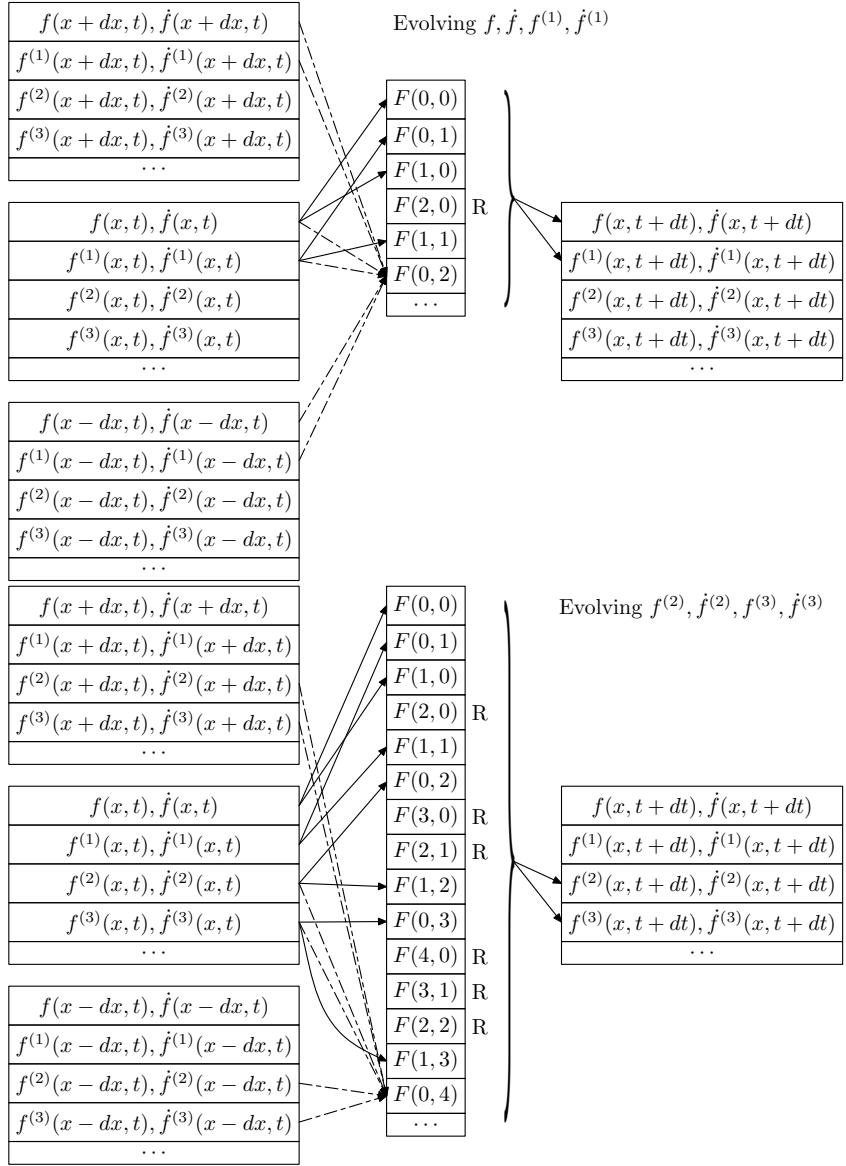

\begin{center}
\leavevmode
\includegraphics[height=0.4\textheight]{diag1.mps}
\includegraphics[height=0.4\textheight]{diag2.mps}
\end{center}
\caption{This diagram illustrates the stacking technique for a second-order-in-time PDE using the pairwise spatial-order evolution. $f^{(n)}$ indicates the $n^{\text{th}}$ spatial derivative of $f$. The solid arrows indicate direct (scaled) assignment. The dot-dashed lines indicate computation using the interpolation procedure detailed in Appendix~\ref{sect:ipoly}. The label ``R'' indicates computation using the recurrence relation derived from the PDE. The right-curly bracket indicates the linear combinations from taking derivatives of the Taylor expansion which defined the recursion coefficients.}
\label{fig:diag}
\end{figure}

An IMDTM scheme where many higher-order coefficients are stored on the grid and interpolation is used only for the highest-order coefficients is numerically unstable. This arises for the same reason that the missing higher-order coefficients cannot be set to zero: a constant error in the order-$n$ coefficient will lead to an error in the order-$(n-1)$ coefficient which grows linearly with time, and so on. In practice, only a couple of orders can be evolved together. You can store only a couple of orders on the grid and use interpolation to generate the coefficients necessary for forward evolution. That can be stable, but lacks the self-consistency constraint. A stable scheme with a self-consistency constraint can be constructed by ``stacking'' a number of such lower-order schemes on top of one another. For example, interpolate using orders zero and one to evolve orders zero and one, interpolate using orders two and three to evolve orders two and three, and so on. This is illustrated in Figure~\ref{fig:diag}. The following examples demonstrate this technique.

\subsection{Example: The Wave Equation}
\label{sect:we}
The salient features of the IMDTM scheme are highlighted with a simple example: the 1+1-dimensional linear wave equation. This example, like the nonlinear equation example which follows, have been chosen because they have known periodic, non-singular, analytic solutions. The implementation uses periodic boundary conditions so that the interpretation of the results is not complicated by details of more-complicated boundary-condition handling. Since each grid point stores an entire tower of derivatives, the boundary points would need to do this as well. How best to use this freedom for other kinds of boundary conditions will be the subject of future research.

For this demonstration, the IMDTM implementation used a grid of $N = 16$ points with the coefficients of the first 14 spatial orders, including the constant term, stored at each point for both the field value and its first time derivative (i.e. $F(0,0)$ through $F(0,13)$ and $F(1,0)$ through $F(1,13)$ are stored at each point). The physical size of the grid is $L = 18$, so $dx = 1.125$ and $dt = 1$ is used for the temporal evolution. In this case, the values of $dx$ and $dt$ are exactly-representable floating-point numbers, although the accuracy is not significantly degraded if this is not the case. The initial conditions were $\phi_0(x) = \cos(2 \pi x / L), \phi_0'(x) = 0$. Unlike traditional explicit temporally-higher-order algorithms (e.g. RK4), IMDTM does not require additional grid copies to hold intermediate grid states\footnote{The method of (spectral) deferred corrections can also be used to produce higher-order temporal evolution for general time-dependent PDEs without storing several intermediate grid copies~\cite{Dutt2000}. As mentioned previously, for systems of conservation laws, ADER schemes can generate higher-order temporal evolution without storing intermediate grid copies.}. Nevertheless, the performance of the IMDTM implementation is compared to an RK4-based scheme. For the particular system chosen here, a second-order stencil (with $dx \approx 0.08)$ is accurate to approximately 4 significant digits, and an eighth-order spatial stencil is accurate to over 13 digits. The $dt$ values for the RK4 scheme were chosen to match the single-core running time to that of the IMDTM implementation\footnote{All code was compiled using GNU g++ version 4.4.4 given the flags: -O3 -march=native -msse -mfpmath=sse -DNDEBUG.}.

\begin{figure}[h!tb]
\begin{center}
\leavevmode
\includegraphics[width=\textwidth]{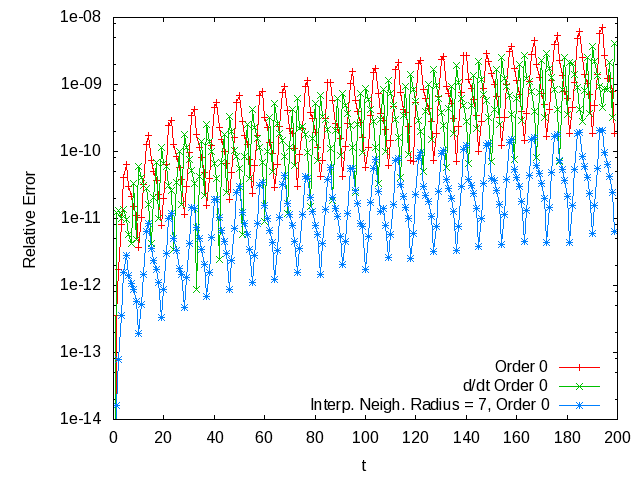}
\end{center}
\caption{The relative error compared to the analytic solution for the system described in Section~\ref{sect:we} with 1 DTM coefficient per point per temporal order. The interpolation uses a six-grid-point radius (13 points per neighborhood). Coefficients of order 13 were the highest-order coefficients used. Also shown is the relative error using an interpolation neighborhood radius of seven (15 points per neighborhood).}
\label{fig:errcmp0}
\end{figure}

Next, consider the most-basic stable use of the IMDTM scheme: only the function values (the zeroth-order coefficients) are stored on the grid. The error\footnote{The average error at each time slice is computed by taking the mean of the base-ten logarithm of the relative error at each point.} of this method in solving the wave equation is illustrated in Figure~\ref{fig:errcmp0}. In such a configuration, the scheme becomes a kind of fancy finite-differencing algorithm. As expected, using a larger interpolation radius leads to a more-accurate evolution; this is exactly equivalent to using a higher-order finite-differencing stencil. 

\begin{figure}[h!tb]
\begin{center}
\leavevmode
\includegraphics[width=\textwidth]{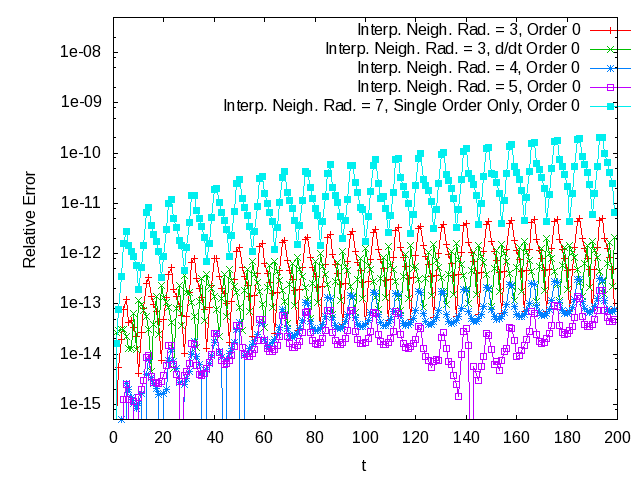}
\end{center}
\caption{The relative error compared to the analytic solution for the system described in Section~\ref{sect:we} with 2 DTM coefficients per point per temporal order. This plot shows cases where the interpolation uses the two highest orders from all points within a three-grid-point radius (7 points per neighborhood), a four-grid-point radius (9 points per neighborhood), and a five-grid-point radius (11 points per neighborhood). Coefficients of order 14 were the highest-order coefficients used. For comparison, also shown is the case where only the zeroth-order components are stored as a seven-grid-point interpolation radius is used. As can be seen, using two orders per grid point yields a much-more-accurate evolution than using only one.}
\label{fig:errcmp1}
\end{figure}

For the case where only the function values are stored, and using a unit radius for interpolation (meaning each point and its two immediate neighbors contribute to the interpolation), it is practical to write down the explicit time-stepping rule. Expanding Equations~\ref{eq:ipoly} and~\ref{eq:a-as-eval} (using Equations~\ref{eq:deriv-using-log} and~\ref{eq:prod-rule} as noted)\footnote{The equations in Appendix~\ref{sect:ipoly} use $\alpha$ to represent the multiset of spatial derivative orders; to apply those formulas to the 1+1-dimensional case here, note that $\alpha$ always has unit cardinality, and its only integer element runs from zero to infinity. In practice, the sum must be truncated at some finite order (determined by the maximum possible order of the interpolating polynomial, which is determined by the number of neighbors and the number of stored spatial derivative orders used for interpolation).}, and using the recurrence relations derived from the wave equation's PDE (given in Equation~\ref{eq:wave-rec}) with Equations~\ref{eq:tsprop00} and~\ref{eq:tsprop10}, yields the following time-stepping rule:
\begin{eqnarray}
	f(t+dt, x) & = & \frac{1}{36 dx^2}\left(-18 \left(dt^2-2 dx^2\right) f(t,x)+ \right . \\ & & \hspace{2em} \left . dt \left(9 dt f(t,x-dx)+9 dt f(t,x+dx)-2 dt^2 \dot{f}(t,x) + \right . \right . \nonumber \\ & & \hspace{2em} \left . \left . 36 dx^2 \dot{f}(t,x)+dt^2 \dot{f}(t,x-dx)+dt^2 \dot{f}(t,x+dx)\right)\right) \nonumber \\
	\dot{f}(t+dt, x) & = & \frac{1}{12 dx^2} \left ( -12 dt f(t,x)+6 dt f(t,x-dx)+6 dt f(t,x+dx)- \right . \\ & & \hspace{2em} \left . 2 dt^2 \dot{f}(t,x)+12 dx^2 \dot{f}(t,x)+dt^2 \dot{f}(t,x-dx)+dt^2 \dot{f}(t,x+dx) \right ) \nonumber
\label{eq:up0}
\end{eqnarray}
where the overdot indicates a time derivative. Note that because the interpolating polynomials are fitting only one value from each of three points, they are only quadratic functions, and so can only estimate the first and second spatial derivatives at each point. So, at each point, $\{ F(0,0), F(1,0) \}$ are given and $\{ F(0,1), F(0,2), F(1,1), F(1,2) \}$ are estimated using the interpolating polynomials\footnote{To be clear, two interpolating polynomials are constructed: One to estimate the spatial derivatives of $f(t,x)$ and one to estimate the spatial derivatives of $\dot{f}(x,t)$.}. From the recurrence relation given by Equation~\ref{eq:wave-rec}, the values of $\{ F(2,0), F(3,0) \}$ can be calculated. The sums in Equations~\ref{eq:tsprop00} and~\ref{eq:tsprop10} can be truncated accordingly. Because the PDE is linear, these same equations are obeyed for each spatial derivative order. The point of this paper, which is that there can be a benefit to storing and evolving spatial derivative values, is demonstrated next. Figure~\ref{fig:errcmp1} shows the case where two orders are stored on the grid and both of those orders are used for the interpolation procedure. Even with a smaller interpolation radius, the result is much more accurate than in the preceding case. As expected, using more coefficients from closer grid points yields a more accurate result.

\begin{figure}[h!tb]
\begin{center}
\leavevmode
\includegraphics[width=\textwidth]{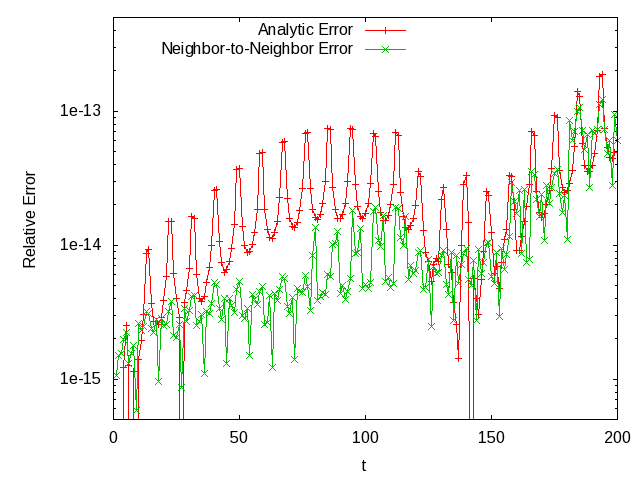}
\end{center}
\caption{The relative error compared to the analytic solution for the system described in Section~\ref{sect:we} with 14 DTM coefficients per point per temporal order. This plot shows cases where the interpolation uses a five-grid-point radius (11 points per neighborhood). Each pair of spatial orders is evolved together, and, because the PDE is linear, independently of all of the other orders. Coefficients of order 25 were the highest-order coefficients used (coefficients of order 14 were the highest-order coefficients used for the first pair, and so on). The neighbor-to-neighbor constraint-violation error, also shown, becomes less than the analytic error.}
\label{fig:errcmp2}
\end{figure}

The IMDTM implementation yields very-high precision evolution while taking large time steps\footnote{This linear case does not form a fair basis for benchmark-like comparison with other numerical evolution schemes, such as an RK4 implementation, in general, because evaluating the recursion relations necessary for the computation of the nonlinear terms can add significant expense, whereas the same is not true for schemes which work only with the function values.}. This could be advantageous for parallel implementations on non-shared-memory clusters, because each time step involves slow synchronizing communication. As for the construction of a scheme with the aforementioned neighbor-to-neighbor self-consistency constraint, the suggested mechanism of ``stacking'' evolution schemes for the coefficients in order pairs is demonstrated in Figure~\ref{fig:errcmp2}. As can be seen by comparing the neighbor-to-neighbor constraint-violation error and the analytic error, the solution is more self-consistent than it is accurate. Comparing the absolute errors instead of the relative errors does not change this conclusion.

\subsection{Example: A Strongly-Nonlinear Equation}
\label{sect:nle}
\begin{figure}[h!tb]
\begin{center}
\leavevmode
\includegraphics[width=\textwidth]{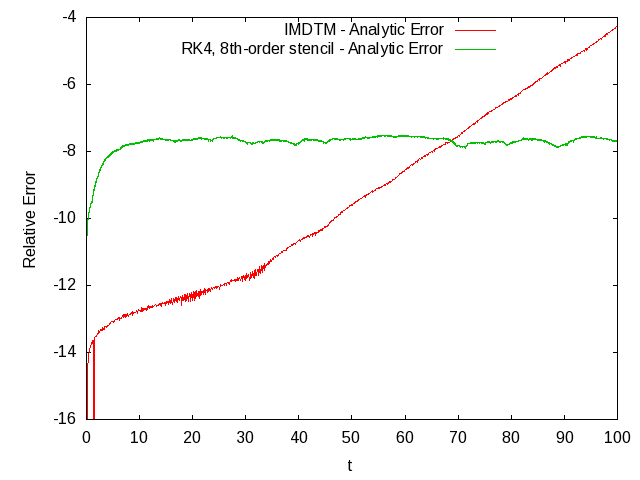}
\end{center}
\caption{The relative error compared to the analytic solution for the system described in Section~\ref{sect:nle} with 2 DTM coefficients per point per temporal order. This plot shows cases where the interpolation uses a five-grid-point radius (11 points per neighborhood). Coefficients of order 17 were the highest-order coefficients used. The RK4 scheme shown used $N = 156$, $dx = 0.28125$ and $dt = 0.001/2650 = 3.77358e-07$, but using $dt = 0.001$ with the RK4 scheme yields an almost-identical curve.}
\label{fig:errcmpnl1}
\end{figure}

While solving the linear wave equation allows the demonstration of several features of the IMDTM scheme, it is not restricted to linear problems. As an example, the IMDTM scheme can be used to evolve the strongly-nonlinear PDE in Equation~\ref{eq:nlpde}. This equation has a non-singular, periodic solution~\cite{He2006}:
\begin{equation}
	-2\sqrt{2}a + \frac{6\sqrt{2}a}{2 \pm \cos(2ax - 8a^3t)} \, .
	\label{eq:nlsltn}
\end{equation}
where $a$ is a free real parameter. Applying the multivariate DTM (see Table~\ref{tab:basemdtm}) to Equation~\ref{eq:nlpde} yields:
\begin{eqnarray}
	H(k, h) & = & \sum_{m=0}^k \sum_{n=0}^h F(k-m, n) F(m, h-n) \\
	G(k, h) & = & \sum_{m=0}^k \sum_{n=0}^h (h-n+1) H(k-m, n) F(m, h-n+1) \\
	F(k, h) & = & -\frac{1}{k} \left (G(k-1, h) + (h+3)(h+2)(h+1) F(k-1, h+3) \right )
\end{eqnarray}
where $H(k,h)$ corresponds to the $f^2$ factor, and $G(k,h)$ corresponds to the $f^2 \frac{\partial f}{\partial x}$ term. Caching the computed values of $H(k,h)$ is essential to an efficient implementation. The stable evolution of Equation~\ref{eq:nlsltn} as the solution of Equation~\ref{eq:nlpde} requires a smaller $dx$ and a smaller $dt$ compared to those used for the linear wave equation. For this demonstration, the IMDTM implementation used a grid of $N = 78$ points. The physical size of the grid is $L = 43.875$, so $dx = 0.56250$ and $dt = 0.001$ is used for the temporal evolution. The initial conditions were $\phi_0(x) = -2\sqrt{2}a + \frac{6\sqrt{2}a}{2 + \cos(2ax)}, a = \pi/L$. As should be expected given the range of the sums in the recurrence relations, unlike for linear PDEs, computing the higher-order terms is significantly more expensive than computing the lower-order terms. Figure~\ref{fig:errcmpnl1} shows how the IMDTM scheme, when only two orders are stored on the grid, compares to an RK4 implementation performing the same amount of computational work and using the same number of grid degrees-of-freedom (not counting the auxiliary grid copies used for the intermediate RK steps). However, reducing the RK4 time step to $dt = 0.001$, identical to that used for the IMDTM scheme, does not significantly degrade the RK4 accuracy. While the IMDTM scheme is less stable than the RK4 scheme, it is far more accurate for tens of thousands of time steps. Although the IMDTM scheme does need extra buffers to efficiently evaluate the nonlinear terms, for this equation it still uses significantly less memory than the RK4 implementation. Although it might seem ``more fair'' to compare to a RK16 scheme, the memory use of an RK16 implementation would be almost 10 times that of the IMDTM scheme, and would likely be unsuitable for practical, large-scale implementation.

\begin{figure}[h!tb]
\begin{center}
\leavevmode
\includegraphics[width=\textwidth]{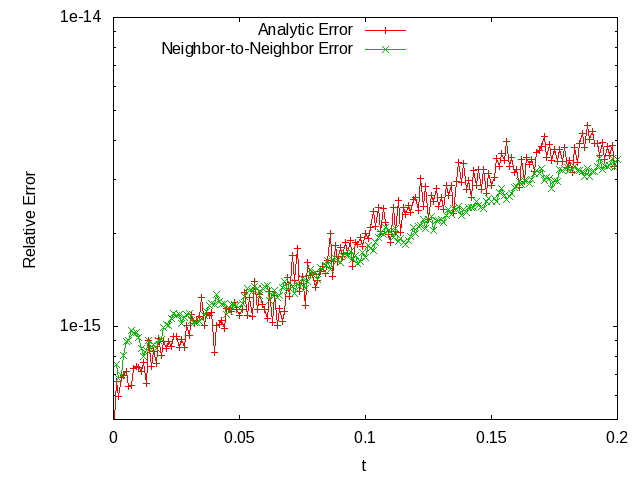}
\end{center}
\caption{The relative error compared to the analytic solution for the system described in Section~\ref{sect:nle} with 17 DTM coefficients per point per temporal order. This plot shows cases where the interpolation uses a five-grid-point radius (11 points per neighborhood). Each pair of spatial orders is evolved together. Coefficients of order 32 were the highest-order coefficients used (coefficients of order 17 were the highest-order coefficients used for the first pair, and so on).}
\label{fig:errcmpnl2}
\end{figure}

\begin{figure}[h!tb]
\begin{center}
\leavevmode
\includegraphics[width=\textwidth]{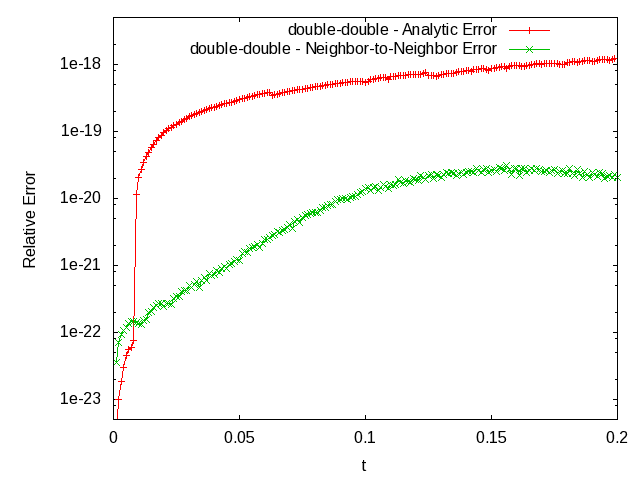}
\end{center}
\caption{The relative error compared to the analytic solution for the system described in Section~\ref{sect:nle} with 17 DTM coefficients per point per temporal order computed using the double-double type from Bailey's qd package~\cite{Hida2001}, which has twice the working precision of an IEEE double. This plot shows cases where the interpolation uses a five-grid-point radius (11 points per neighborhood). Each pair of spatial orders is evolved together. Coefficients of order 32 were the highest-order coefficients used (coefficients of order 17 were the highest-order coefficients used for the first pair, and so on).}
\label{fig:errcmpnl3}
\end{figure}

Pair-wise order evolution (as illustrated in Figure~\ref{fig:diag}) was used in order to generate a scheme with a neighbor-to-neighbor self-consistency constraint, and the results were similar to those for the linear case (shown in Figure~\ref{fig:errcmpnl2}). The growth of the analytic error over the neighbor-to-neighbor constraint-violation error is exacerbated when computing at higher precision. Specifically, when using the double-double type from Bailey's qd package~\cite{Hida2001}, which has twice the precision of an IEEE double, after some hundreds of time steps, the analytic error becomes greater than the neighbor-to-neighbor constraint-violation error. This is demonstrated in Figure~\ref{fig:errcmpnl3}.

\section{Future Work}
The work presented here leaves a lot of room for future research. The following are some of the items which deserve further investigation:
\begin{itemize}
\item A procedure for performing a stability analysis of an IMDTM scheme should be developed.
\item The IMDTM should apply naturally to problems with nontrivial boundary conditions, so long as the multiple power series expansion can be generated in a self-consistent manner. The exact conditions under which this is possible should be investigated.
\item DTM can be applied to boundary-value problems by leaving some of the lower-order coefficients free and then solving for them in a way consistent with the boundary constraints after expressions for the higher-order coefficients have been obtained~\cite{Hassan2007}~\cite{Hassan2009}. Similarly, it may be possible to extend IMDTM to handle multivariate boundary-value problems and evolution problems with constraints.
\item Alternate methods for computing the power series at time $t + \delta t$ from those at time $t$ which are more stable than the method presented here might exist. For example, it has been demonstrated that using Pad\'e approximants can increase the stability of an iterated DTM scheme~\cite{Haq2010}~\cite{Peker2010}. It may also be possible to formulate the IMDTM as an implicit scheme.
\item A reusable framework for developing IMDTM codes is under development, and its interface definitions and features will be presented in a future publication.
\end{itemize}

\section{Conclusion}
In this paper, the IMDTM PDE evolution scheme for initial-value problems has been introduced. Methods to stabilize the evolution and efficiently use polynomial interpolation to enable iterative forward propagation have been detailed. It is an explicit evolution scheme, and it tends to be able to take larger time steps compared to other explicit schemes (e.g. RK4) while performing a similar amount of computational work. For many equations, the IMDTM scheme will also use less memory than the corresponding RK-like scheme. The larger time steps should allow for substantially-reduced communication on non-shared-memory machines. Because the algorithm is local and explicit, it should scale strongly to the largest conceivable supercomputers. Furthermore, IMDTM can naturally be used to construct a verifiably-self-consistent code in the sense that it has an internal self-consistency constraint which can be used to measure the quality of the solution. Not only does this often allow the user to know when the solution can no longer be trusted, it can allow the IMDTM solution to be used as a baseline for verifying other codes. More generally, this paper demonstrates that efficient numerical schemes can be constructed which store and evolve higher-order spatial derivative information, and that doing so can confer significant practical advantages.

\section*{Acknowledgments}
Most of this work was done while I was supported by the United States Department of Energy Computational Science Graduate Fellowship, provided under grant DE-FG02-97ER25308. I am now supported under grant DE-AC02-06CH11357. I thank Richard Easther for his advice, feedback, and for reading an early draft of this paper. I also thank Christopher Gilbreth, Jim Stewart, Bill Rider and Matt Norman for providing corrections and useful suggestions.

\bibliographystyle{plain}
\bibliography{imdtm}

\vfill
\begin{center}
\scriptsize
\framebox{\parbox{3.2in}{The submitted manuscript has been created by UChicago Argonne, LLC, Operator of Argonne National Laboratory ("Argonne").  Argonne, a U.S. Department of Energy Office of Science laboratory, is operated under Contract No. DE-AC02-06CH11357.  The U.S. Government retains for itself, and others acting on its behalf, a paid-up, nonexclusive, irrevocable worldwide license in said article to reproduce, prepare derivative works, distribute copies to the public, and perform publicly and display publicly, by or on behalf of the Government.}} \normalsize \end{center}

\appendix
\section{Interpolating Polynomial}
\label{sect:ipoly}
The interpolating polynomial can be derived from the multipoint Taylor expansion as defined by L\'opez and Temme~\cite{Lopez2004}. Let $\{ x_i^{(1)}, x_i^{(2)}, \ldots, x_i^{(m)} \}$ be the set of $m$ unique coordinate values in dimension $i \in \{ 1, \ldots, D \} $. In one dimension, a function $f(x)$ has the multipoint expansion:
\begin{equation}
	f(x) = \sum_{n = 0}^\infty \left ( \sum_{j=1}^m \frac{\prod_{k=1,k \neq j}^m (x - x^{(k)})}{\prod_{k=1,k \neq j}^m (x^{(j)} - x^{(k)})} a_{n,j} \right ) \prod_{k=1}^m (x - x^{(k)})^n
\end{equation}
where
\begin{eqnarray}
	a_{n,j} = \frac{1}{n!} \frac{d^n}{dx^n} \left [ \frac{f(x)}{\prod_{s=1,s \neq j}^m (x - x^{(s)})^n} \right ]_{x = x^{(j)}} \nonumber \\ + \sum_{k=1,k \neq j}^m \frac{1}{(n-1)!} \frac{d^{n-1}}{dx^{n-1}} \left [ \frac{f(x)}{(x - x^{(j)}) \prod_{s=1,s \neq k}^m (x - x^{(s)})^n} \right ]_{x = x^{(k)}}
\end{eqnarray}
which can be written as:
\begin{equation}
	a_{n,j} = \sum_{k=1}^m \frac{1}{(n - (k \neq j))!} \frac{d^{n - (k \neq j)}}{dx^{n - (k \neq j)}} \left [ \frac{f(x)}{(x - x^{(j)})_{k \neq j} \prod_{s=1,s \neq k}^m (x - x^{(s)})^n} \right ]_{x = x^{(k)}} \, .
\end{equation}
This can be generalized in the usual way to a multidimensional function by recursive expansion, yielding:
\begin{eqnarray}
	f(x) = \sum_{\alpha} \left ( \sum_{\beta \in \{ 1, \ldots, m \}^D} \frac{\prod_{i=1}^D \prod_{k=1, k \neq \beta_i}^m (x_i - x_i^{(k)})}{\prod_{i=1}^D \prod_{k=1, k \neq \beta_i}^m (x_i^{(\beta_i)} - x_i^{(k)})} a_{\alpha,\beta} \right ) \nonumber \\ \times \prod_{i=1}^D \prod_{k=1}^m (x_i - x_i^{(k)})^{\alpha_i} \, ,
\end{eqnarray}
Expanding the final, non-constant binomial term gives:
\begin{eqnarray}
	f(x) = \sum_{\alpha} \left ( \sum_{\beta \in \{ 1, \ldots, m \}^D} \frac{\prod_{i=1}^D \sum_{\gamma \in \{0, \ldots, m\} \times \{0,1\}^{m-1}, \origbar\gamma\origbar=m-1} x_i^{\gamma_0} \prod_{k=1, k \neq \beta_i}^m (-x_i^{(k)})^{\gamma_k}}{\prod_{i=1}^D \prod_{k=1, k \neq \beta_i}^m (x_i^{(\beta_i)} - x_i^{(k)})} a_{\alpha,\beta} \right ) \nonumber \\ \times \prod_{i=1}^D \prod_{k=1}^m \left ( \sum_{p=0}^{\alpha_i} {\alpha_i \choose p} (-x_i^{(k)})^{\alpha_i-p} x_i^p \right ) \, .
\end{eqnarray}
Note that the product of two sequences of length $N$ is their convolution:
\begin{equation}
	\left ( \sum_{n=0}^N a_n x^n \right ) \left ( \sum_{n=0}^N b_n x^n \right ) = \sum_{i=0}^{2N} \left ( \sum_{i=0}^n a_i b_{n-i} \right ) x^n
\end{equation}
and the same applies for three sequences:
\begin{eqnarray}
	\left ( \sum_{n=0}^N a_n x^n \right ) \left ( \sum_{n=0}^N b_n x^n \right ) \left ( \sum_{n=0}^N c_n x^n \right ) = \left ( \sum_{i=0}^{2N} \left ( \sum_{i=0}^n a_i b_{n-i} \right ) x^n \right ) \left ( \sum_{n=0}^N c_n x^n \right ) \nonumber \\ = \sum_{i=n}^{3N} \left ( \sum_{j=0}^n  \left ( \sum_{i=0}^j a_i b_{j-i} \right ) c_{n-j} \right ) x^n \, ,
\end{eqnarray}
from which the pattern can be seen.
\begin{eqnarray}
	f(x) = \sum_{\alpha} \left ( \sum_{\beta \in \{ 1, \ldots, m \}^D} \frac{\prod_{i=1}^D \sum_{\gamma \in \{0, \ldots, m\} \times \{0,1\}^{m-1}, \origbar\gamma\origbar=m-1} x_i^{\gamma_0} \prod_{k=1, k \neq \beta_i}^m (-x_i^{(k)})^{\gamma_k}}{\prod_{i=1}^D \prod_{k=1, k \neq \beta_i}^m (x_i^{(\beta_i)} - x_i^{(k)})} a_{\alpha,\beta} \right ) \nonumber \\ \times \prod_{i=1}^D \left ( \sum_{p=0}^{m\alpha_i} x_i^p \sum_{q_1=0}^p \sum_{q_2=0}^{q_1} \cdots \sum_{q_m=0}^{q_{m-1}} {\alpha_i \choose q_m} (-x_i^{(m)})^{\alpha_i-q_m} \right . \nonumber \\ \left . \vphantom{\sum_{p=0}^{m\alpha_i}} \times {\alpha_i \choose {q_{m-1} - q_m}} (-x_i^{(m-1)})^{\alpha_i-(q_{m-1}-q_m)} \cdots {\alpha_i \choose {p - q_1}} (-x_i^{(1)})^{\alpha_i-(p-q_1)} \right ) 
\end{eqnarray}
The contribution to the DTM coefficient $F(h)$ is then easily derived from the $x^h$ coefficient of $f(x)$:
\begin{eqnarray}
	f(x) = \sum_\alpha \sum_{\beta \in \{ 1, \ldots, m \}^D} \frac{a_{\alpha,\beta}}{\prod_{i=1}^D \prod_{k=1, k \neq \beta_i}^m (x_i^{(\beta_i)} - x_i^{(k)})} \nonumber \\ \times \prod_{i=1}^D \sum_{p=0}^{m^2\alpha_i} x_i^p \sum_{j=0}^p \left ( \sum_{\gamma \in \{0,1\}^{m-1}, \origbar\gamma\origbar=m-1-(p-j)} \prod_{k=1, k \neq \beta_i}^m (-x_i^{(k)})^{\gamma_k} \right ) \nonumber \\ \times \left ( \sum_{q_1=0}^j \sum_{q_2=0}^{q_1} \cdots \sum_{q_m=0}^{q_{m-1}} {\alpha_i \choose q_m} (-x_i^{(m)})^{\alpha_i-q_m} \right . \nonumber \\ \left . \vphantom{\sum_{p=0}^{m\alpha_i}} \times {\alpha_i \choose {q_{m-1} - q_m}} (-x_i^{(m-1)})^{\alpha_i-(q_{m-1}-q_m)} \cdots {\alpha_i \choose {j - q_1}} (-x_i^{(1)})^{\alpha_i-(j-q_1)} \right )
	\label{eq:ipoly}
\end{eqnarray}
and
\begin{equation}
	a_{\tilde{\alpha},\beta} = \sum_{\gamma \in \{ 1, \ldots, m \}^D} \frac{1}{\alpha!} \frac{\partial^{\origbar\alpha\origbar}}{\partial x^\alpha} \left [ \frac{f(x)}{\prod_{i=1}^D (x_i - x_i^{(\beta_i)})_{\gamma_i \neq \beta_i} \prod_{k=1,k \neq \gamma_i}^m (x_i - x_i^{(k)})^{\tilde{\alpha}_i}} \right ]_{x = x^{(\gamma)}}
\end{equation}
where $\alpha = \tilde{\alpha} - \{\gamma_i \neq \beta_i\}$. The product rule for a general partial derivative is~\cite{Hardy2006}:
\begin{equation}
{\partial^{k_1+\cdots+k_n} \over \partial x_1^{k_1}\,\cdots\,
\partial x_n^{k_n}} (uv) 
=\sum_{\ell_1=0}^{k_1}\cdots\sum_{\ell_n=0}^{k_n}
{k_1 \choose \ell_1}\cdots {k_n \choose \ell_n}
{\partial^{\ell_1+\cdots+\ell_n} u
\over \partial x_1^{\ell_1}\,\cdots\,\partial x_n^{\ell_n}}
\cdot {\partial^{k_1-\ell_1+\cdots+k_n-\ell_n} v
\over \partial x_1^{k_1-\ell_1}\,\cdots
\,\partial x_n^{k_n-\ell_n}} \, .
\end{equation}
So we can expand the expression for $a_{\tilde{\alpha},\beta}$:
\begin{eqnarray}
	a_{\tilde{\alpha},\beta} = \sum_{\gamma \in \{ 1, \ldots, m \}^D} \frac{1}{\alpha!} \sum_{\delta, \forall i\,\delta_i \leq \alpha_i} {\alpha \choose \delta} \frac{\partial^{\origbar\delta\origbar}}{\partial x^\delta} [ f(x) ]_{x = x^{(\gamma)}} \nonumber \\ \times \frac{\partial^{\origbar\alpha - \delta\origbar}}{\partial x^{\alpha - \delta}} \left [ \frac{1}{\prod_{i=1}^D (x_i - x_i^{(\beta_i)})_{\gamma_i \neq \beta_i} \prod_{k=1,k \neq \gamma_i}^m (x_i - x_i^{(k)})^{\tilde{\alpha_i}}} \right ]_{x = x^{(\gamma)}} \, ,
\end{eqnarray}
and the factor for each dimension can be separated:
\begin{eqnarray}
	a_{\tilde{\alpha},\beta} = \sum_{\gamma \in \{ 1, \ldots, m \}^D} \frac{1}{\alpha!} \sum_{\delta, \forall i\,\delta_i \leq \alpha_i} {\alpha \choose \delta} \frac{\partial^{\origbar\delta\origbar}}{\partial x^\delta} [ f(x) ]_{x = x^{(\gamma)}} \nonumber \\ \times \prod_{i=1}^D \frac{\partial^{\origbar\alpha_i - \delta_i\origbar}}{\partial x_i^{\alpha_i - \delta_i}} \left [ \frac{1}{(x_i - x_i^{(\beta_i)})_{\gamma_i \neq \beta_i} \prod_{k=1,k \neq \gamma_i}^m (x_i - x_i^{(k)})^{\tilde{\alpha_i}}} \right ]_{x = x^{(\gamma)}} \, .
	\label{eq:a-as-eval}
\end{eqnarray}
Because $\frac{d}{dx} \ln f = \frac{1}{f} \frac{df}{dx}$ we can write $\frac{df}{dx} = f \frac{d}{dx} \ln f$, and so:
\begin{eqnarray}
	\frac{\partial}{\partial x_i} \left [ \frac{1}{(x_i - x_i^{(\beta_i)})_{\gamma_i \neq \beta_i} \prod_{k=1,k \neq \gamma_i}^m (x_i - x_i^{(k)})^{\tilde{\alpha}_i}} \right ]_{x = x^{(\gamma)}} = \nonumber \\ \left [ \frac{1}{(x_i - x_i^{(\beta_i)})_{\gamma_i \neq \beta_i} \prod_{k=1,k \neq \gamma_i}^m (x_i - x_i^{(k)})^{\tilde{\alpha}_i}} \left ( \left ( -\frac{1}{x_i - x_i^{(\beta_i)}} \right )_{\gamma_i \neq \beta_i} - \tilde{\alpha}_i \sum_{k=1, k \neq \gamma_i}^m \frac{1}{x_i - x_i^{(k)}} \right ) \right ]_{x = x^{(\gamma)}} \, .
	\label{eq:deriv-using-log}
\end{eqnarray}
The second factor in Equation~\ref{eq:a-as-eval} can be evaluated recursively using Equation~\ref{eq:deriv-using-log} and the product rule:
\begin{equation}
	\frac{d^n}{dx^n} uv = \sum_{k=0}^n {n \choose k} \frac{d^{n-k}u}{dx^{n-k}} \frac{d^kv}{dx^k} \, .
	\label{eq:prod-rule}
\end{equation}

\section{Additional Multivariate DTM Recurrences}
\label{sect:mdtmrec}
The following table shows how to apply the Multidimensional DTM to some common nonlinear terms. Once the formulas for the basic arithmetic operations are known, deriving the formulas for functions that are solutions of an ODE is straightforward~\cite{Forbes1986}~\cite{Kedem1980} (just apply the DTM transformation rules to the defining ODE). A large table of these recurrences can be found in an appendix of the Fast Forward Automated Differentiation Library (FFADLib) User Manual~\cite{Tsukanov2000}.

\begin{center}
\begin{longtable}{|l|p{0.75\textwidth}|}
\hline
\bf{Original Function} & \bf{Transformed Function} \\
\hline
$w(x) = y(x)/z(x)$ & $W(k) = \frac{1}{Z(0)} \left [ Y(k) - \sum_{l_1=0}^{k_1} \cdots \sum_{l_n=0}^{k_n} Z(l) W(k-l) \right ] $ \\
$w(x) = \sqrt{y(x)}$ & $W(k) = \frac{1}{2W(0)} \left [ Y(k) - \sum_{l_1=0}^{k_1} \cdots \sum_{l_n=0}^{k_n} W(l) W(k-l) \right ]$ \\
$w(x) = \exp(y(x))$ & $W(k) = \frac{1}{k_a} \sum_{l_1=0}^{k_1} \cdots \sum_{l_a=0}^{k_n-1} \cdots \sum_{l_n=0}^{k_n}$ \newline \hspace*{2em} $(k_a - l_a) W(l) Y(k-l)$ \\
$w(x) = \ln(y(x))$ & $W(k) = \frac{1}{Y(0)} \left [ Y(k) - \frac{1}{k_a} \sum_{l_1=0}^{k_1} \cdots \sum_{l_a=0}^{k_n-1} \cdots \sum_{l_n=0}^{k_n} \right .$ \newline \hspace*{2em} $\left . \vphantom{\sum_{l_1=0}^{k_1}} (k_a - l_a) Y(l) W(k-l) \right ]$ \\
$w(x) = y(x)^s$ & $W(k) = \frac{1}{Y(0)} \left [ s W(0) Y(k) + \frac{1}{k_a} \sum_{l_1=0}^{k_1} \cdots \sum_{l_a=0}^{k_n-1} \cdots \sum_{l_n=0}^{k_n} \right .$ \newline \hspace*{2em} $\left . \vphantom{\sum_{l_1=0}^{k_1}} (k_a - l_a) [ s W(l) Y(k-l) - Y(l) W(k-l) ] \right ]$ \\
$w(x) = \sin(y(x))$ & $W(k) = \frac{1}{k_a} \sum_{l_1=0}^{k_1} \cdots \sum_{l_a=0}^{k_n-1} \cdots \sum_{l_n=0}^{k_n}$ \newline \hspace*{2em} $(k_a - l_a) W_{\cos{}}(l) Y(k-l)$ \\
$w(x) = \cos(y(x))$ & $W(k) = - \frac{1}{k_a} \sum_{l_1=0}^{k_1} \cdots \sum_{l_a=0}^{k_n-1} \cdots \sum_{l_n=0}^{k_n}$ \newline \hspace*{2em} $(k_a - l_a) W_{\sin{}}(l) Y(k-l)$ \\
\hline
\end{longtable}
\end{center}

For the self-referential recurrences, note that $W(0)$ is always determined by applying the original function to the zeroth-order coefficients of the respective power series. In addition, any terms which would cause $W(k)$ to depend on itself should be omitted; their appearance is only for notational convenience. For some of the entries it is necessary to pick an index, $a$, for which the sum only goes from $0$ to $k_a-1$. The choice of this index is arbitrary, although it generally must be such that $k_a \neq 0$.
\end{document}